\documentclass[12pt,draft]{extartic}

\usepackage{amsmath,amsthm,amssymb,calc}
\usepackage[dvips]{graphics, graphicx}

\usepackage[english]{babel}

\def\eps{\varepsilon}

\newcounter{num}[section]

\newcommand{\Th}{\refstepcounter{num}
{\bf Theorem \arabic{section}.\arabic{num} }}

\newcommand{\Lemma}{\refstepcounter{num}
{\bf Lemma \arabic{section}.\arabic{num} }}

\newcommand{\Pred}{\refstepcounter{num}
{\bf Proposition \arabic{section}.\arabic{num} }}

\newcommand{\Cor}{\refstepcounter{num}
{\bf Corollary \arabic{section}.\arabic{num} }}

\newcommand{\Note}{\refstepcounter{num}
{\it Note \arabic{section}.\arabic{num} }}

\newcommand{\St}{\refstepcounter{num}
{\bf Statement \arabic{section}.\arabic{num} }}

\newcommand{\Def}{\refstepcounter{num}
{\it Definition \arabic{section}.\arabic{num} }}

\newcommand{\Proof}{{\bf Proof. }}

\def\eps{\varepsilon}
\def\_phi{\varphi}

\def\a{\alpha}
\def\d{\delta}
\def\l{\lambda}

\def\L{\Lambda}
\def\m{\times}
\def\t{\tilde}

\def\ov{\overline}
\def\z{{\mathbb Z}}

\def\Z_N{{\mathbb Z}_N}
\def\Z{{\mathbb Z}}

\def\Span{{\rm Span\,}}

\author{Shkredov I.D.}
\title{On sets with small doubling
\footnote{This work was supported by RFFI grant no.
06-01-00383, President's of Russian Federation grant N 1726.2006.1
and INTAS (grant no. 03--51--5-70).}}
\date{}
\begin{document}
\maketitle

\begin{center}
 Annotation.
\end{center}

{\it \small Let $G$ be an arbitrary Abelian group
and let $A$ be a finite subset of $G$.
$A$ has small additive doubling if $|A+A| \le K|A|$ for some $K>0$.
These sets were studied in papers of G.A. Freiman, Y. Bilu, I. Ruzsa, M.C.--Chang, B. Green and
T.Tao.
In the article we prove that if we have some minor restrictions on $K$
then for any set with small doubling there exists a set $\L$, $\L \ll_\eps K \log |A|$
such that $|A\cap \L| \gg |A| / K^{1/2 + \eps}$, where $\eps > 0$.
In contrast to the previous results our theorem is nontrivial
for large $K$.
For example one can take $K$ equals $|A|^\eta$,
where
$\eta>0$.
We use an elementary method in our proof.
}
\\
\\
\\

\refstepcounter{section}
\label{introduction}

{\bf \arabic{section}. Introduction.}

Let $G$ be
an arbitrary Abelian group with additive
group operation $+$.
Suppose that $A,B$ are two finite subsets of $G$
and define their sumset $A+B$
to be the set of all pairwise sums $a+b$ with $a\in A$, $b \in B$.
Let $\log$ stand for the logarithm to base $2$.

Suppose that $A$ is a set such that $|A+A| \le K |A|$,
where $K \ge 1$ is small
(for example $K = \log \log |A|$ or $K=2$).
These sets are called sets with {\it small doubling}.
The properties of such sets were studied in papers
\cite{Freiman,Bilu,Ruzsa_Freiman,Ch_Fr,Ruzsa_Freiman_groups,Ruzsa_Green_Fr,Green_pol,Green_Sanders,Green_Tao_U3,Green_Tao_flat}.
G.A. Freiman (see \cite{Freiman}) proved the following wonderful result on
the structure of these sets.

Recall that a set $Q \subseteq G$ is called a {\it $d$--dimensional arithmetic progression} if
 $$
    Q = \{ n_0 + n_1 \l_1 + \dots + n_d \l_d ~:~ 0 \le \l_i < m_i \} \,,
 $$
 where $m_i, n_i \in \Z$ and $m_i \ge 0$.

 Let $G=\Z$.

 \Th {\bf (Freiman)}
 \label{t:Freiman}
 {\it
    Let $K\ge 1$ be a real number, and
    $A \subseteq \Z$ be a finite set.
    Let also $|A+A| \le K|A|$.
    Then there exist numbers $d = d(K)$ and $C = C(K)$ depend on $K$ only
    and $d$--dimensional arithmetic progression $Q$ such that $|Q| \le C |A|$
    and $A\subseteq Q$.
 }

 The functions $d = d(K)$ and $C = C(K)$ were studied in \cite{Ruzsa_Freiman,Ch_Fr}.
 In paper \cite{Ch_Fr} M.--\,C. Chang proved that
 $d = O (K^2 \log^2 K)$ and $C = \exp ( O(K^2 \log^2 K) )$
 (as usual we use $X = O(Y)$ or $X \ll Y$
  to denote an estimate of the form $X \le M Y$
  for some absolute constant $M$).

 Let $n$ be a positive integer.
 Sets with small doubling in groups $G = \left( \z / 2\z \right)^n$
 were considered in
 \cite{Ruzsa_Freiman_groups,Sanders_Fr,Green_pol,Green_Sanders,Green_Tao_flat}.
 For example we formulate a theorem from \cite{Ruzsa_Freiman_groups}.
 Note that $\left( \z / 2\z \right)^n$ is a vector space.

\Th
\label{t:Ruzsa10}
{\it
    Let $K \ge 1$ be a real number.
    Let $A \subseteq \left( \z / 2\z \right)^n$ be a set such that $|A+A| \le K |A|$.
    Then $A$ is contained in a subspace $H$
    with $|H| \le K^2 2^{K^4} |A|$.
}

    There are another structural results on sets with small doubling.
    Let $A$ be a set with small doubling and $d$ is a small positive integer.
    Is it true that $A$ has large {\it intersection} with
    some $d$---dimensional arithmetic progression?
    It is known that there is a positive answer at the question
    and
    we give two examples of such results.
    In \cite{Green_Sanders} the following theorem was proved.

\Th
\label{t:Green_Sanders}
{\it
    Let $K \ge 1$ be a real number.
    Suppose that $A \subseteq \left( \z / 2\z \right)^n$ is a set
    such that $|A+A| \le K |A|$.
    Then there exists a subspace $H$
    such that $|H| \ll K^{O(1)} |A|$
    and
    $|A\cap H| \gg \exp( - K^{O(1)} ) |A|$.
}

    Finally, in recent paper \cite{Green_Tao_flat} B. Green and T. Tao
    proved the following theorem.

\Th
\label{t:Green_Tao_flat}
{\it
    Let $K \ge 1$ be a real number.
    Let $A \subseteq \left( \z / 2\z \right)^n$ be a set such that
    $|A+A| \le K |A|$.
    Then there exists a subspace $H$ and $x \in \left( \z / 2\z \right)^n$
    such that $|H| \gg K^{-O(\sqrt{K})} |A|$
 and
    $|A\cap (x+H)| \ge \frac{1}{2K} |H|$.
}

    Let us formulate our main result.

    Let $E = \{ e_1, \dots, e_{|E|} \} \subseteq G$ be a finite set.
    By $\Span E$ denote the set
    $\Span E = \{~ \sum_{i=1}^{|E|} \eps_i e_i ~:~ \eps_i \in \{ -1,0,1 \}~ \}$.

\Th
\label{t:main}
{\it
    Let $G$ be an Abelian group.
    Let $K,\eps$ be real numbers, $\eps \in (0,1/2]$,
    $A \subseteq G$ be a finite set, $|A| \ge 2^{32/\eps}$,
    $1 \le K \le \min \{~ ( 2^{-58} \eps^{-4} \frac{|A|}{\log |A|} )^{(3/2 + \eps)^{-1}}, |A|^{\eps} ~\}$.
    Let also $A$ contains at least $|A|^3 / K$
    quadruples with
    $a_1 + a_2 = a_3 + a_4$.
    Then there exists a set $\L$ such that
    $|\Span \L \cap A| \ge \frac{1}{2} \cdot \frac{|A|}{K^{1/2 + \eps}}$
    and
    $|\L| \le 2^{30} \eps^{-2} K \log |A|$.
}

    It is easy to see that the number $K$ in Theorems \ref{t:Freiman},\ref{t:Ruzsa10},\ref{t:Green_Sanders},\ref{t:Green_Tao_flat}
    cannot be too large.
    For example Theorem \ref{t:Green_Tao_flat} is trivial if
    $K \gg \left( \frac{\log |A|}{\log \log |A|} \right)^2$.
    In contrast to these results our Theorem \ref{t:main}
    is nontrivial for large $K$
    (for example one can take $K = |A|^{\eta}$, where $\eta > 0$ is a sufficiently small number).
    On the other hand the cardinality of the set $\L$ depends on $|A|$.
    This fact  differences our main result from Theorems
    \ref{t:Freiman}---\ref{t:Green_Tao_flat}.

    This paper is organized as follows.
    In \S \ref{connected_sets} we study so--called
    "connected sets"\, in Abelian groups.
    We prove that any such set has large intersection with $\Span \L$
    for some small set $\L$
    (in more detail see Propositions \ref{pr:small_basis}, \ref{pr:almost_basis}).
    Besides in the section we show that any set contains large connected subset.
    These two facts imply Theorem \ref{t:main}.
    We give its proof in \S \ref{proof_main}.
    In the last section we discuss some relations between our definition
    of connectedness and a graph--theoretical definition of connectedness
    from paper \cite{Ruzsa_Elekes}.

The author is grateful to Professor N.G. Moshchevitin
for constant attention to this work.

\refstepcounter{section}
\label{connected_sets}

{\bf \arabic{section}. On connected sets in Abelian groups.}

Let $G$ be an arbitrary Abelian group with additive
group operation $+$.
Let $A\subseteq G$ be a finite set, and $k\ge 2$ be a positive integer.
By $T_k (A)$ denote the following number
$$
    T_k (A) := | \{ a_1 + \dots + a_k = a'_1 + \dots + a'_k  ~:~ a_1, \dots, a_k, a'_1,\dots,a'_k \in A \} | \,.
$$
Denote by the same letter $A$ the characteristic function of $A$.
Thus $A(x)=1$ if $x\in A$ and $A(x)=0$ otherwise.
We shall write $\sum_x$ instead of $\sum_{x\in G}$.

\Def
\label{b_connectedness}
    Let $k\ge 2$ be a positive integer, and $\beta\in [0,1]$ be a real number.
    Suppose that $A$ is a finite nonempty set $A\subseteq G$.
    $A$ is called {\it $\beta$--connected of degree $k$} if
    there is an absolute constant $C\in (0,1]$
    such that for any $B\subseteq A$, $|B| \ge \beta |A|$ we have
    \begin{equation}\label{ineq:b_connectedness}
        T_k (B) \ge C^{2k} \left( \frac{|B|}{|A|} \right)^{2k} T_k (A) \,.
    \end{equation}
    If $\beta = 0$ then $A$ is {\it connected of degree $k$}.

    The class of connected sets is wide enough.
    On the one hand very structured sets like arithmetic progressions,
    multidimensional arithmetic progressions, subspaces are connected sets
    (see Corollary \ref{cor:subspaces-a_conn} below).
    On the other hand any so--called {\it dissociated} set
    (see Definition \ref{def:diss})
    belongs to this class.
    Other examples of connected of degree $k$ sets will be
    considered in section \ref{other_connectedness}.


\Def
\label{def:convolution}
    Let $f,g : G \to \mathbb{R}$ be arbitrary functions.
    Denote by $(f*g) (x)$ the function
\begin{equation}\label{f:*-conv}
    (f*g) (x) = \sum_s f(s) g(x-s) \,.
\end{equation}
    Clearly, $(f*g) (x) = (g*f) (x)$, $x\in G$.
    By $(f \circ g) (x)$ denote the function
\begin{equation}\label{f:circ-concv}
    (f \circ g) (x) = \sum_s f(s) g(s-x) \,.
\end{equation}
    Obviously, $(f \circ g) (x) = (g \circ f) (-x)$, $x\in G$.

Suppose that  $A,B \subseteq G$ are any sets.
Then $(A*B) (x) \neq 0$ iff $x\in A+B$
and $(A \circ B) (x) \neq 0$ iff  $x\in A-B$.
Hence $T_2 (A) = \sum_x (A*A)^2 (x)$.
Further denote by $*_{k}$ the composition of $k$ operations $*$, $k\ge 1$.
Then $T_k (A) = \sum_x (A *_{k-1} A)^2 (x)$, $k\ge 2$.
Since
$$
    T_2 (A)
        :=
            |\{ a_1 + a_2 = a'_1 + a'_2  ~:~ a_1, a_2, a'_1, a'_2 \in A \}|
                =
                    |\{ a_1 - a'_1 = a'_2 - a_2 ~:~ a_1, a_2, a'_1, a'_2 \in A \}|
$$
it follows that
$T_2 (A) = \sum_x (A \circ A)^2 (x)$.

    Let $f: G \to \mathbb{R}$ be a function.
    By $T_k (f)$ denote the quantity
    $T_k (f) = \sum_x (f *_{k-1} f)^2 (x)$.
    Let us prove the following simple lemma.

\Lemma
\label{l:T_k(B,A)}
{\it
    Let $p_1$,$p_2$ be positive integers, and $k_1=2^{p_1}$, $k_2=2^{p_2}$.
    Let also
    $f_1, \dots, f_{k_1}, g_1, \dots, g_{k_2} : G \to \mathbb{R}$ be functions.
    Then
    $$
        \left| \sum_x (f_1 * \dots * f_{k_1}) (x) \cdot (g_1 * \dots * g_{k_2}) (x) \right|
            \le
    $$
    \begin{equation}\label{f:T_k(B,A)}
            \le
                (T_{k_1} (f_1))^{1/2k_1} \dots (T_{k_1} (f_{k_1}))^{1/2k_1}
                    (T_{k_2} (g_1))^{1/2k_2} \dots (T_{k_2} (g_{k_2}))^{1/2k_2} \,.
    \end{equation}
}
\Proof
    First of all let us suppose that $k_1 = k_2 = k = 2^p$, where $p$ is a positive integer.
    We prove the lemma by induction.
    Put
    $\sigma = \sum_x (f_1 * \dots * f_k) (x) \cdot (g_1 * \dots g_k) (x)$.
    Using the Cauchy--Schwartz inequality, we get
    \begin{equation}\label{f:begin_l}
        \sigma^2 \le \sum_x (f_1 * \dots * f_k)^2 (x) \cdot \sum_x (g_1 * \dots * g_k)^2 (x) = \sigma_1 \sigma_2 \,.
    \end{equation}
    Consider the sum $\sigma_1$.
    By definitions of $*,\circ$, we obtain
    $$
        \sigma_1 = \sum_x ( (f_1 \circ f_1) * \dots * (f_{2^{p-1}} \circ f_{2^{p-1}})) (x)
                                \cdot
                                    ( (f_{2^{p-1}+1} \circ f_{2^{p-1}+1}) * \dots * (f_k \circ f_k)) (x)
    $$
    By the induction hypothesis, we have
    \begin{equation}\label{}
        \sigma_1 \le (T_{2^{p-1}} (f_1 \circ f_1))^{1/k} \dots (T_{2^{p-1}} (f_k \circ f_k))^{1/k} \,.
    \end{equation}
    Besides,
    $T_{2^{p-1}} (f_1 \circ f_1) = T_k (f_1)$.
    Hence
    \begin{equation}\label{tmp:23:05_18}
        \sigma_1 \le (T_k (f_1))^{1/k} \dots (T_k (f_k))^{1/k} \,.
    \end{equation}
    In the same way
    \begin{equation}\label{tmp:23:05_18+}
        \sigma_2 \le (T_k (g_1))^{1/k} \dots (T_k (g_k))^{1/k} \,.
    \end{equation}
    Combining (\ref{tmp:23:05_18}), (\ref{tmp:23:05_18+}) and
    (\ref{f:begin_l}), we obtain that (\ref{f:T_k(B,A)})
    holds.

    Let now $k_1=2^{p_1}$, $k_2=2^{p_2}$, and $p_1 \neq p_2$.
    Put
    $\sigma' = \sum_x (f_1 * \dots * f_{k_1}) (x) \cdot (g_1 * \dots g_{k_2}) (x)$.
    Using the Cauchy--Schwartz inequality, we get
    \begin{equation}\label{f:begin_l'}
        \sigma'^2 \le \sum_x (f_1 * \dots * f_{k_1})^2 (x) \cdot \sum_x (g_1 * \dots * g_{k_2})^2 (x) = \sigma'_1 \sigma'_2 \,.
    \end{equation}
    Using (\ref{f:T_k(B,A)}) for $\sigma'_1$, $\sigma'_2$, we have
    \begin{equation}\label{}
        |\sigma'| \le (T_{k_1} (f_1))^{1/2k_1} \dots (T_{k_1} (f_{k_1}))^{1/2k_1}
                           (T_{k_2} (g_1))^{1/2k_2} \dots (T_{k_2} (g_{k_2}))^{1/2k_2} \,.
    \end{equation}
    This completes the proof.

    Let us derive a corollary from Lemma \ref{l:T_k(B,A)}.
    Let $n$ be positive integer, $q$ be a prime and let $G$ be $(\Z /q\Z)^n$.
    As was noted above $G$ is a vector space.

\Cor
\label{cor:subspaces-a_conn}
{\it
    Let $n$, $p$ be positive integers, $k=2^p$, $q$ be a prime, and $G = (\Z /q\Z)^n$.
    Let also $P$ be a subspace of $G$.
    Then $P$ is a connected of degree $k$ set and (\ref{ineq:b_connectedness})
    is true for $C=1$.
}
\\
\Proof
    Let $B \subseteq P$ be a set
    and let $\sigma (B) := \sum_x (B * P *_{k-2} P) (x) \cdot (P *_{k-1} P) (x)$.
    The sum $\sigma (B)$ equals the number of solutions
    of the equation
    $b + p_2 + \dots + p_k = p'_1 + \dots + p'_k$, where
    $b\in B$ and $p_2, \dots, p_k, p'_1, \dots, p'_k \in P$.
    Since $P$ is a subspace of $(\Z /q\Z)^n$ it follows that
    $b + p_2 + \dots + p_k - p'_1 - \dots - p'_k \in P$.
    Hence $\sigma (B) \ge |B| |P|^{2k-2}$.
    In particular $\sigma (P) = T_k (B) \ge |P|^{2k-1}$.
    Since $T_k (P) \le |P|^{2k-1}$ it follows that $T_k (P) = |P|^{2k-1}$.
    Using Lemma \ref{l:T_k(B,A)} with $f_1 = B$,
    $f_2 = \dots = f_k = g_1 = \dots = g_k = A$, we get
    \begin{equation}\label{}
        \sigma^{2k} (B) \le T_k (B) \cdot T^{2k-1}_k (A) \,.
    \end{equation}
    Combining the last inequality and the lower bound for $\sigma (B)$, we obtain
    $$
        T_k (B)
            \ge
                \frac{\sigma^{2k} (B)}{T^{2k-1}_k (A)}
                    \ge
                        \frac{|B|^{2k} |P|^{(2k-2)2k} }{|P|^{(2k-1)^2}}
                            =
                                \left( \frac{|B|}{|P|} \right)^{2k} |P|^{2k-1}
                                    =
                                        \left( \frac{|B|}{|P|} \right)^{2k} T_k (P) \,.
    $$
    This completes the proof.

    Thus very structured sets like subspaces are connected sets.
    Consider another examples of connected sets.

    We need in the following definition (see \cite{Rudin_book} or \cite{Ch_Fr}).

\Def
\label{def:diss}
We say that $\L = \{ \l_1, \dots, \l_{|\L|} \} \subseteq G$ is
{\it dissociated} if the equality
   \begin{equation}\label{f:def_diss}
        \sum_{i=1}^{|\L|} \eps_i \l_i = 0 \,,
   \end{equation}
where $\eps_i \in \{ -1,0,1 \}$ implies that all $\eps_i$ are equal to zero.

If $\L$ is a dissociated set then there exists a good upper bound for $T_k (\L)$
(see \cite{Rudin_book} and also \cite{Green_Chang1,Rudin}).

\St
\label{st:Rudin_T_k}
{\it
    There is an absolute constant $M>0$ such that
    for any dissociated set $\L \subseteq G$ and any positive integer $k\ge 2$, we have
    \begin{equation}\label{tmp:3_2}
        T_k (\L) \le M^{k} k^k |\L|^k \,,
    \end{equation}
    where $M \le 288$.
}

Any connected  of degree $k$ set has the following property.

\Pred
\label{pr:small_basis}
{\it
    Let $k\ge 2$ be a positive integer.
    Suppose that $A\subseteq G$ is a connected  of degree $k$ set
    and for $C>0$ inequality (\ref{ineq:b_connectedness})
    holds.
    Then there exists a set $\L \subseteq A$,
    $
        |\L| \le 288 C^{-2} k \frac{|A|^2}{T^{1/k}_k (A)}
    $
    such that any $a\in A$
    can be expressed in the form
    \begin{equation}\label{f:r=log}
        a = \sum_{i=1}^{|\L|} \eps_i \l_i \,,
    \end{equation}
    where $\eps_i \in \{ -1,0,1 \}$.
}
\\
\Proof
   Let $\L$ be a maximal dissociated subset of $A$.
   Prove that any $a\in A$ can be expressed in the form
   \begin{equation}\label{tmp:4:47}
        a = \sum_{i=1}^{|\L|} \eps_i \l_i \,,
   \end{equation}
   where $\eps_i \in \{ -1,0,1 \}$.
   If $a=0$ then (\ref{tmp:4:47}) holds.
   Let $a$ be an arbitrary element of $A \setminus \L$, $a\neq 0$.
   Consider all equations  $\sum_{i=1}^{|\L|+1} \eps_i \t{\l}_i = 0$, where
   $\t{\l}_i \in \L \bigsqcup \{ a \}$ and
   $\eps_i \in \{ -1,0,1 \}$, $i\in \{ 1,2, \dots, |\L|+1 \}$.
   If all these equations are trivial, i.e. we have
   $\eps_i = 0$, $i\in \{ 1,2, \dots, |\L|+1 \}$ then
   we obtain a contradiction with the maximality of $\L$.
   It follows that there exists non--trivial equation
   $\eps a + \sum_{i=1}^{|\L|} \eps_i \l_i = 0$, $\eps, \eps_i \in \{ -1,0,1 \}$
   such that not all
   $\eps, \eps_i$ are equal to zero.
   Note that $\eps\neq 0$.
   Whence any  $a\in A$ is involved in some equation (\ref{f:r=log}).

   Let us prove that $|\L| \le 288 C^{-2} k \frac{|A|^2}{T^{1/k}_k (A)}$.
   Using Statement \ref{st:Rudin_T_k},
   we have
   $T_k (\L) \le (288)^k k^k |\L|^k$.
   On the other hand, the set $A$ is connected of degree $k$.
   Hence
   $T_k (\L) \ge C^{2k} (|\L| / |A| )^{2k} \cdot T_k (A)$.
   It follows that $|\L| \le 288 C^{-2} k \frac{|A|^2}{T^{1/k}_k (A)}$.
   This completes the proof.

We need in a more delicate definition of connectedness.

\Def
\label{b_connectedness+}
    Let $k\ge 2$ be a positive integer,
    and $\beta_1, \beta_2 \in [0,1]$ be real numbers, $\beta_1 \le \beta_2$.
    Nonempty set $A\subseteq G$ is called
    {\it $(\beta_1, \beta_2)$--connected of degree $k$} if
    there exists an absolute constant $C\in (0,1]$
    such that for any $B\subseteq A$,
    $\beta_1 |A| \le |B| \le \beta_2 |A|$ we have
    \begin{equation}\label{ineq:b_connectedness+}
        T_k (B) \ge C^{2k} \left( \frac{|B|}{|A|} \right)^{2k} T_k (A) \,.
    \end{equation}

Clearly, any $\beta$--connected of degree $k$ set
is a $(\beta, \beta_2)$--connected of degree $k$,
where $\beta_2 \in [\beta,1]$ is an arbitrary number.
Nevertheless we have the following weak analog of
Proposition \ref{pr:small_basis}
for $(\beta_1, \beta_2)$--connected of degree $k$ sets.

\Pred
\label{pr:almost_basis}
{\it
    Let $k\ge 2$ be a positive integer, $0<\beta_1 \le \beta_2$ be real numbers.
    Let also $A\subseteq G$ be a $(\beta_1, \beta_2)$--connected of degree $k$ set
    and
    for $C>0$ inequality (\ref{ineq:b_connectedness})
    holds.
    Suppose that $\beta_2 \ge \beta_1 + 1/|A|$,
    $T_k (A) \ge 2^{14k} C^{-2k} k^k |A|^k$
    and
    $|A| \ge 1/\beta_1$.
    Then there exists a set $\L \subseteq A$
    such that
    \begin{equation}\label{est:L+}
        |\L| \le 2^{13} C^{-2} k \frac{|A|^2}{T^{1/k}_k (A)} \,,
    \end{equation}
    and $|\Span \L \cap A| \ge (1-\beta_1) |A|$.
}
\\
\Proof
   The proof of the proposition is a sort of inductive process.
   Let $\L_1$ be a dissociated subset of $A$ such that
   $|\Span \L_1 \cap A| \ge (1-\beta_1) |A|$.
   Clearly, there exists such set $\L_1$,
   for example
   one put $\L_1$ to be a maximal dissociated subset of $A$.
   Let $l = 2^{13} C^{-2} k \frac{|A|^2}{T^{1/k}_k (A)}$.
   If $|\L_1| \le l$ then the proposition is proved.
   Suppose that $|\L_1| > l$.
   Let $\L'_1 \subseteq \L_1$ be an arbitrary set of the cardinality $l$.
   Obviously, that $\L'_1$ is a dissociated set.
   Consider the set $A_1 = A\setminus \L'_1$.
   If $|A_1| < (1-\beta_1) |A|$ then we stop our algorithm.
   If $|A_1| \ge (1-\beta_1) |A|$
   then
   let $\L_2$ be a dissociated subset of $A_1$ such that
   $|\Span \L_2 \cap A_1| \ge (1-\beta_1) |A|$.
   Suppose that $|\L_2| \le l$.
   Then $|\Span \L_2 \cap A| \ge |\Span \L_2 \cap A_1| \ge (1-\beta_1) |A|$
   and we are done.
   It follows that $|\L_2| > l$.
   Let $\L'_2 \subseteq \L_2$ be an arbitrary set of the cardinality $l$
   and consider the set $A_2 = A_1 \setminus \L'_2$.
   An so on.
   We get the sets $A_0 = A, A_1, A_2, \dots, A_s$
   and disjoint dissociated sets
   $\L'_1, \dots, \L'_s$ from  $A$.
   We have $|A_s| < (1-\beta_1) |A|$.
   Since for all $l=1,2,\dots,s$ the following holds $A_l = A \setminus \bigsqcup_{i=1}^l \L'_i$
   it follows that
   $\sum_{i=1}^s |\L'_i| = |A| - |A_s| > \beta_1 |A|$.
   Let $B = \bigsqcup_{i=1}^s \L'_i$.
   Then $|B| > \beta_1 |A|$.
   We can remove some elements from $\L'_s$ and
   assume that the cardinality of $\bigsqcup_{i=1}^s \L'_i$
   equals $[\beta_1 |A|] + 1$.
   Denote by the same letter $B$ our modified set.
   We have $B\subseteq A$ and $|B| \ge \beta_1 |A|$.
   Since $\beta_2 \ge \beta_1 + 1/|A|$ it follows that $|B| \le \beta_2 |A|$.
   By assumption the set $A$ is $(\beta_1, \beta_2)$ connected of degree $k$.
   Hence
   \begin{equation}\label{Rom_3}
        T_k (B) \ge C^{2k} \beta^{2k}_1 T_k (A) \,.
   \end{equation}
   On the other hand
   \begin{equation}\label{Rom_3'}
        T_k (B)
            \le
                T_k \left( \bigsqcup_{i=1}^s \L'_i \right)
                    =
                        \sum_{i_1,\dots,i_k = 1}^s ~  \sum_{j_1,\dots,j_k = 1}^s
                            \sum_x (\L'_{i_1} * \dots * \L'_{i_k}) (x) \cdot (\L'_{j_1} * \dots * \L'_{j_k}) (x) \,.
   \end{equation}
   Using Lemma \ref{l:T_k(B,A)}, Statement \ref{st:Rudin_T_k}
   and (\ref{Rom_3'}), we get
   \begin{equation}\label{Rom_4}
        T_k (B) \le s^{2k} \max_{i=1,\dots,s} T_k (\L'_i) \le s^{2k} (288)^k k^k l^k \,.
   \end{equation}
   By assumption $T_k (A) \ge 2^{14k} C^{-2k} k^k |A|^k$.
   Whence $|A| \ge 2^{14} C^{-2} k \frac{|A|^2}{T^{1/k}_k (A)} = 2l$
   and $s\ge 2$.
   Since $\bigsqcup_{i=1}^{s-1} \L'_i \subseteq B$ and $|A| \ge 1/\beta_1$
   it follows that $sl /2 \le (s-1) l \le |B| \le 2 \beta_1 |A|$.
   Hence $s \le 4 \beta_1 |A| / l$.
   Combining the last inequality and (\ref{Rom_4}), we have
   $$
        T_k (B) \le 2^{4k} \beta^{2k}_1 (288)^k k^k \frac{|A|^{2k}}{l^k} \,.
   $$
   This contradicts with (\ref{Rom_3})
   and we obtain the required result.

Let us prove now that any $A\subseteq G$ contains
some large $(\beta_1, \beta_2)$--connected of degree $k$ set.
We begin with some notation.

\Def
\label{def:zeta}
    Let $A \subseteq G$ be an arbitrary finite set, $|A| \ge 2$,
    and $k\ge 2$ be a positive integer.
    By $\zeta_k (A)$ denote the quantity
    $$
        \zeta_k = \zeta_k (A) := \frac{\log T_k (A)}{\log |A|} \,.
    $$

In other words $T_k (A) = |A|^{\zeta_k}$.
Clearly, for any set $A$, we have $k\le \zeta_k (A) \le 2k-1$.

    Let $A \subseteq G$ be a finite set, $|A| = m \ge 2$,
    $p$ be a positive integer, and $k=2^p$.
    Write $\zeta$ for $\zeta_k (A)$.

\Th
\label{t:connected_subset}
{\it
    Let $\beta_1, \beta_2 \in (0,1)$ be real numbers, $\beta_1 \le \beta_2$.
    Then there exists a set $A' \subseteq A$ such that \\
    $1)~$ $A'$ is $(\beta_1, \beta_2)$--connected of degree $k$ set
            such that (\ref{ineq:b_connectedness+}) holds for any $C \le 1/32$.\\
    $2)~$ $|A'| \ge m \cdot 2^{\frac{\log ( (2k-1)/ \zeta)}{\log (1+\kappa)} \log (1-\beta_2)}$,
                    where $\kappa = \frac{\log ((1-\beta_1)^{-1})}{\log m} (1-16 C)$.\\
    $3)~$ $\zeta_k (A') \ge \zeta_k (A)$.
}
\\
\Proof
    Let $C\le 1/32$ be a real number.
    The proof of Theorem \ref{t:connected_subset} is a sort algorithm.
    If $A$ is $(\beta_1, \beta_2)$--connected of degree $k$
    and (\ref{ineq:b_connectedness+}) is true with the constant $C$
    then there is nothing to prove.
    Suppose that $A$ is not
    $(\beta_1, \beta_2)$--connected of degree $k$ set (with the constant $C$).
    Then there exists a set $B\subseteq A$, $\beta_1 |A| \le |B| \le \beta_2 |A|$
    such that (\ref{ineq:b_connectedness+}) does not hold.
    Note that $|A|>2$.
    Let $\ov{B} = A\setminus B$.
    We have
    $$
        T_k (A) = \sum_x (A *_{k-1} A)^2 (x)
            =
    $$
    \begin{equation}\label{f:t_1.1}
            =
                \sum_x (B * A *_{k-2} A) (x) (A *_{k-1} A) (x)
                    +
                        \sum_x (\ov{B} * A *_{k-2} A) (x) (A *_{k-1} A) (x)
                            =
                                \sigma_1 + \sigma_2 \,.
    \end{equation}
    Using Lemma \ref{l:T_k(B,A)} with $f_1 = B$, $f_2 = \dots = f_k = g_1 = \dots = g_k = A$, we obtain
    \begin{equation}\label{f:t_1.tmp1}
        \sigma_1^{2k} \le T_k (B) \cdot T^{2k-1}_k (A) \,.
    \end{equation}
    In the same way
    \begin{equation}\label{f:t_1.tmp2}
        \sigma_2^{2k} \le T_k (\ov{B}) \cdot T^{2k-1}_k (A) \,.
    \end{equation}
    Let $c_B = |B|/ |A|$.
    Combining $T_k (B) < C^{2k} c^{2k}_B T_k (A)$ and (\ref{f:t_1.tmp1}),
    we have $\sigma_1 < C c_B T_k (A)$.
    Using the last inequality, (\ref{f:t_1.1}) and (\ref{f:t_1.tmp2}), we get
    \begin{equation}\label{f:t_1.2}
        T_k (\ov{B}) > T_k (A) ( 1 - C c_B )^{2k} \,.
    \end{equation}
    Let $\ov{\zeta} = \zeta_k (\ov{B})$, $b = |B|$ and $\ov{b} = |\ov{B} | = m - b$.
    Using (\ref{f:t_1.2}), we obtain
    $$
        \ov{\zeta} \log \ov{b} > \zeta \log m + 2k \log (1-C c_B) \,.
    $$
    Hence
    $$
        \ov{\zeta}
            >
                \frac{\zeta \log m + 2k \log (1-C c_B)}{\log \ov{b}}
                    =
                         \frac{\zeta \log m + 2k \log (1-C c_B)}{\log m + \log (1-b/m)}
                            =
                                \frac{\zeta + 2k \frac{\log (1-C c_B)}{\log m}}{1 + \frac{\log (1-c_B)}{\log m}}
                                    \ge
    $$
    $$
                                    \ge
                                        \left( \zeta + 2k \frac{\log (1-C c_B)}{\log m} \right)
                                            \left( 1 - \frac{\log (1-c_B)}{\log m} \right)
                                                \ge
                                                    \zeta + \zeta \frac{\log ((1-c_B)^{-1})}{\log m} (1-16 C)
                                                        \ge
    $$
    \begin{equation}\label{f:t_1.3}
                                                        \ge \zeta (1+\frac{\log ((1-\beta_1)^{-1})}{\log m} (1-16 C))
                                                            = \zeta (1+\kappa) \,,
    \end{equation}
    where $\kappa = \frac{\log ((1-\beta_1)^{-1})}{\log m} (1-16 C) > 0$.
    Besides, by the definition of  $(\beta_1, \beta_2)$--connectedness of degree $k$,
    we have
    \begin{equation}\label{f:t_1.4}
        |\ov{B}| \ge (1-\beta_2) m = (1-\beta_2) |A| \,.
    \end{equation}
    Thus if the set $A$ is not $(\beta_1, \beta_2)$--connected of degree $k$
    then there is a set $\ov{B} \subseteq A$ such that
    (\ref{f:t_1.3}), (\ref{f:t_1.4}) hold.
    Put $A_1 = \ov{B}$ and apply the arguments above to $A_1$.
    And so on.
    We get the sets $A_0 = A, A_1, A_2, \dots, A_s$.
    Clearly, for any $A_i$, we have $\zeta (A_i) \le 2k-1$.
    Using this and (\ref{f:t_1.3}), we obtain that the total number of steps
    of our algorithm
    does not exceed
    $\frac{\log ( (2k-1)/ \zeta)}{\log (1+\kappa)}$.
    At the last step of the algorithm, we find the set $A'=A_s \subseteq A$
    such that $A'$ is $(\beta_1, \beta_2)$--connected of degree $k$
    and $\zeta_k (A') \ge \zeta = \zeta_k (A)$.
    Thus $A'$ has the properties $1)$ and $3)$
    of the theorem.
    Let us prove $2)$.
    Using (\ref{f:t_1.4}), we obtain
    $$
        |A'| \ge (1-\beta_2)^s m \ge m \cdot 2^{\frac{\log ( (2k-1)/ \zeta)}{\log (1+\kappa)} \log (1-\beta_2)} \,.
    $$
    This completes the proof.

\Cor
\label{cor:dense_sets}
{\it
    Let $G$ be an Abelian group, $\eps$, $\d$ be real numbers,
    $\eps \in (0,1/8]$, $\d \in (0,1]$,
    $\d \ge |G|^{-\eps}$  and
    let $A \subseteq G$ be a set, $|A| \ge \d |G| \ge 2$.
    Let also $p$ be a positive integer, $k=2^p$,
    and
    $\beta_1, \beta_2 \in (0,1)$ be real numbers, $\beta_1 \le \beta_2$, $\beta_1 \le 1 - |A|^{-2\eps}$.
    Then there exists a set $A' \subseteq A$ such that \\
    $1)~$ $A'$ is $(\beta_1, \beta_2)$--connected of degree $k$
            and (\ref{ineq:b_connectedness+}) is true for any
            $C \le 1/32$.\\
    $2)~$ $|A'| \ge |G| \cdot \d^{ \left(\frac{2}{2k-1} + 32 \eps \right) \cdot \frac{\log (1-\beta_2)}{\log (1-\beta_1)} + 1}$.\\
    $3)~$ $\zeta_k (A') \ge \zeta_k (A)$.\\
    In particular, if $\beta_2 = \beta_1$, $k=2$ and $\eps = 1/8$
    then the cardinality of $|A'|$ is at least $\d^{6} |G|$.
}
\\
\Proof
Using Theorem \ref{t:connected_subset} with $C=1/32$, we find a set
$A' \subseteq A$ with properties $1)$---$3)$ guaranteed by the theorem.
Let us prove that
$|A'| \ge |G| \cdot \d^{\left(\frac{2}{2k-1} + 32 \eps \right)\cdot \frac{\log (1-\beta_2)}{\log (1-\beta_1)} + 1}$.
Let $N = |G|$, $m=|A|$, and $\zeta = \zeta_k (A)$.
Clearly, $T_k (A) \ge \d^{2k} N^{2k-1}$.
Hence
\begin{equation}\label{f:tmp:22:45_25}
    \zeta \ge 2k-1 + (2k-1) \frac{\log (1/\d)}{\log N} - \frac{2k}{1-\eps} \frac{\log (1/\d)}{\log N} \,.
\end{equation}
Since $\d \ge N^{-\eps}$ it follows that
\begin{equation}\label{tmp:tmp}
    2k-1 - \zeta \le \frac{\log (1/\d)}{\log N} (1 + 4 k \eps) \quad \mbox{ and } \quad \zeta \ge (2k-1) (1 - 5 \eps) \,.
\end{equation}
By Theorem \ref{t:connected_subset}, we have
$|A'| \ge m \cdot 2^{\frac{\log ( (2k-1)/ \zeta)}{\log (1+\kappa)} \log (1-\beta_2)}$,
where $\kappa = \frac{\log ((1-\beta_1)^{-1})}{2\log m}$.
Using the last inequality, $\beta_1 \le 1 - |A|^{-2\eps}$ and
(\ref{f:tmp:22:45_25}), (\ref{tmp:tmp}), we get
$$
    |A'| \ge m \cdot 2^{ \left(\frac{2}{2k-1} + 32 \eps \right) \cdot
                                \frac{\log (1-\beta_2)}{\log (1-\beta_1)} \cdot \frac{\log \d}{\log N} \log m }
         \ge N \cdot \d^{\left(\frac{2}{2k-1} + 32 \eps \right) \cdot
                                \frac{\log (1-\beta_2)}{\log (1-\beta_1)} + 1} \,.
$$
This completes the proof.

\Note
Certainly, the constant $32$
at the second point of Corollary \ref{cor:dense_sets}
can be decreased.
The constant $2$ in the numerator
of $\frac{2}{2k-1}$ depends on an upper bound for $C$.
If $C$ is less than $1/32$ then the number $2$
is also decreases.

\refstepcounter{section}
\label{proof_main}

{\bf \arabic{section}. The proof of main result.}

\Lemma
\label{l:T_k_&_T_2}
{\it
    Let $A$ be a finite nonempty set, and $k$ be a positive integer, $k\ge 2$.
    Then $T_k (A) \ge T^{k-1}_2 (A) / |A|^{k-2}$.
}


\Proof
The proof is trivial.

{\bf The proof of Theorem \ref{t:main}}
Let $m=|A|$, $\beta_1 = 1/2$, $\beta_2 = \beta_1 + 1/\log m$,
$C = \eps 2^{-7}$, $k=2^p$, $p = [\log \ln m] + 1$.
Clearly, $C \le 1/32$.
Using Theorem \ref{t:connected_subset}, we find $A' \subseteq A$
such that $1)$ --- $3)$ hold.
By assumption $T_2 (A) \ge |A|^3 / K$.
Using Lemma \ref{l:T_k_&_T_2}, we get
$T_k (A) \ge T^{k-1}_2 (A) / |A|^{k-2} \ge |A|^{2k-1} / K^{k-1}$.
Whence
\begin{equation}\label{tmp:Rim_1}
    \zeta = \zeta_k (A) \ge 2k-1 - (k-1) \frac{\log K}{\log m} \,.
\end{equation}
Using $K \le m^{\eps}$ and (\ref{tmp:Rim_1}), we obtain
\begin{equation}\label{tmp:Rim_2}
    \zeta \ge (2k-1) \left( 1 - \frac{k-1}{2k-1} \frac{\log K}{\log m} \right)
            \ge (2k-1) \left( 1 - \frac{\eps}{2} \right) \,.
\end{equation}
By $2)$ of Theorem \ref{t:connected_subset}, we have
\begin{equation}\label{tmp:Rim_3}
    |A'| \ge m \cdot 2^{-\frac{\log ( (2k-1)/ \zeta)}{\log (1+\kappa)} \log ((1-\beta_2)^{-1}) } = m 2^{-\sigma} \,,
\end{equation}
where $\kappa = \frac{\log ((1-\beta_1)^{-1})}{\log m} (1-16 C)$.
Let us obtain an upper bound on $\sigma$.
Using simple inequalities $\log (1+x) \le \frac{x}{\ln 2}$, $\log (1+x) \ge \frac{1}{\ln 2} (x - x^2/2)$, $x\ge 0$,
and inequalities $m \ge 2^{32/\eps}$, (\ref{tmp:Rim_1}), (\ref{tmp:Rim_2}), we get
$$
    \sigma
        \le
            \log \left( 1+ \frac{2k-1-\zeta}{\zeta} \right) \frac{\ln 2}{\kappa} (1+\kappa) \log ((1-\beta_2)^{-1})
                \le
$$
$$
                \le
                    \frac{k-1}{2k-1} \cdot \frac{\log K}{\log m} \left( 1 + \eps \right)
                    \frac{\log m}{1-16 C} \left( 1 + \frac{1}{\log m} \right)
                    \frac{\log ((1-\beta_2)^{-1})}{\log ((1-\beta_1)^{-1})}
                        \le
$$
$$
                        \le
                            \log K^{1/2} \left( 1 + \eps \right) (1+32 C) \left( 1 + \frac{8}{\log m} \right)
                                \le \log K^{1/2 + \eps} \,.
$$
Hence $|A'| \ge \frac{m}{K^{1/2 + \eps}}$.
Since $\zeta_k (A') \ge \zeta_k (A)$ it follows that
\begin{equation}\label{tmp:Rim_5}
    T_k (A') \ge \frac{|A'|^{2k-1}}{K^{k-1}} \ge \frac{|A'|^{2k-1}}{K^k} \,.
\end{equation}
Using $|A'| \ge \frac{m}{K^{1/2 + \eps}}$,
(\ref{tmp:Rim_5})
and
$K \le ( 2^{-58} \eps^{-4} \frac{|A|}{\log |A|} )^{(3/2 + \eps)^{-1}}$
it is easy to see that $T_k (A') \ge 2^{14k} C^{-2k} k^k |A'|^k$.
Using Proposition \ref{pr:almost_basis}, we find a set $\L$ such that
$|\Span \L \cap A'| \ge |A'| /2$
and
\begin{equation}\label{tmp:Rim_3'}
    |\L| \le 2^{27} \eps^{-2} k \frac{|A'|^2}{T^{1/k}_k (A')} \,.
\end{equation}
We have
\begin{equation}\label{tmp:Rim_4}
    |\Span \L \cap A| \ge |\Span \L \cap A'| \ge \frac{|A'|}{2} \ge \frac{1}{2} \cdot \frac{m}{K^{1/2 + \eps}} \,.
\end{equation}
Let us prove that $|\L| \le 2^{30} \eps^{-2} K \log m$.
Combining (\ref{tmp:Rim_3'}) and (\ref{tmp:Rim_5}), we get
$$
    |\L| \le 2^{27} \eps^{-2} K k |A'|^{1/k} \le 2^{27} \eps^{-2} K k m^{1/k} \,.
$$
Recall that $k=2^p$, $p = [\log \ln m] + 1$, we finally obtain
$|\L| \le 2^{30} \eps^{-2} K \log m$.
This completes the proof.

\Cor
\label{cor:small_doubling}
{\it
    Let $G$ be an Abelian group .
    Let $K,\eps$ be real numbers, $\eps \in (0,1/2]$,
    $A \subseteq G$ be an arbitrary set, $|A| \ge 2^{32/\eps}$,
    $1 \le K \le \min \{~ ( 2^{-58} \eps^{-4} \frac{|A|}{\log |A|} )^{(3/2 + \eps)^{-1}}, |A|^{\eps} ~\}$.
    Let also $|A+A| \le K |A|$.
    Then there exists a set $\L$ such that
    $|\Span \L \cap A| \ge \frac{1}{2} \cdot \frac{|A|}{K^{1/2 + \eps}}$
    and
    $|\L| \le 2^{30} \eps^{-2} K \log |A|$.
}
\\
\Proof
We have $|A+A| \le K |A|$.
By the Cauchy--Schwartz inequality
$$
    |A|^4 = \left( \sum_x (A * A) (x) \right)^2 \le \sum_x (A * A)^2 (x) \cdot |A+A| \le T_2 (A) \cdot K |A| \,.
$$
Hence $T_2 (A) \ge |A|^3 / K$.
Using Theorem \ref{t:main}, we obtain the required result.

\refstepcounter{section}
\label{other_connectedness}

{\bf \arabic{section}. Another definitions of connectedness.}

    In the section we discuss some relations between our definition
    of connectedness and a graph--theoretical definition of connectedness.

Suppose that $\Gamma = (V,f)$ is a graph,
where $V$ is the set of vertices of $\Gamma$
and $f$ is the characteristic function of
a symmetric subset of $V \m V$.
Let $X,Y \subseteq V$ be arbitrary sets.
By $e(X,Y)$ denote the number of vertices between $X$ and $Y$.
In other words $e(X,Y) = \sum_{x\in X} \sum_{y\in Y} f(x,y)$.
Recall that $\Gamma$ is {\it connected}
if for any vertex $x$, we have $e(x,V\setminus \{ x \}) > 0$.
In \cite{Ruzsa_Elekes} I. Ruzsa and G. Elekes
gave the following definition.

\Def
\label{def:Ruzsa_Elekes}
{
    Let $\a \in (0,1]$ be a real number.
    A graph $\Gamma = (V,f)$ is called
    {\it $\a$--dense--connected} if
    for any partition of the set of vertices into two
    disjoint parts, say $E\bigsqcup F = V$, we have
    $$
        e(E,F) \ge \a |E| |F| \,.
    $$
}

We give an analog of the definition above for subsets of Abelian groups.
Let $G$ be an Abelian group, and
$A\subseteq G$ be a finite set.
In papers \cite{Gow_4,Gow_m,Ch_Fr} the graph of "popular differences"\, of $A$
was considered.
This graph $\Gamma_A = ( V_A, f_A)$
played a significant role in various problems of combinatorial number theory
(see articles \cite{Gow_4,Gow_m,Ch_Fr} and book \cite{Tao_Vu_book}).
The vertex set $V_A$ of the graph $\Gamma_A$ is $A$,
and the function $f_A$ is the characteristic function of the symmetric set of
"popular differences"
\begin{displaymath}
  f (x,y) =
  \left\{ \begin{array}{ll}
    1,                 &            \mbox{ if } \quad |\, \{ x-y = a_1 - a_2 ~:~ a_1,a_2 \in A \} \,| \ge h \,,
                              \\
    0,                 &            \mbox{ otherwise. }
  \end{array} \right.
\end{displaymath}
    Here $h$ is a number, $0 \le h \le |A|$.
    In many problems of combinatorial number theory
    $h$ was taken
    approximately
    $T_2 (A) / |A|^{2}$.
    Thus the function $f(x,y)$ equals $1$ if
    $(A \circ A) (x-y) \ge h$
    and equals $0$ otherwise.
    Ruzsa and Elekes applied $\a$--dense--connected subgraphs of $\Gamma_A$
    to prove some results on sumsets
    (see details in \cite{Ruzsa_Elekes}).

    In the article we define a new (generalized) graph $\Gamma'_A = (V'_A, f'_A)$,
    where $f'_A$ is a symmetric function but not the characteristic function
    of some subset of $V'_A \m V'_A$.
    Put $V'_A := A$ and $f'_A (x,y) := (A\circ A) (x-y)$.
    The constructed graph $\Gamma'_A$ is an "approximation"\, of the graph $\Gamma_A$
    in the sense that the function $f_A$ is a normalized and truncated version of
    the function $f'_A$ : $f_A (x,y) = \theta (f'_A (x,y) / h)$,
    where $\theta$ is a shifted Heaviside's function :~
    $\theta (x) = 1$ if $x\ge 1$ and $\theta (x) = 0$ if $x< 1$.
    Then the graph $\Gamma'_A$ is (generalized) $\a$--dense--connected
    if for any partition of the set of vertices into two
    disjoint parts $E$ and $V$, $E\bigsqcup F = A$, we have
\begin{equation}\label{f:pre_strong}
    e(E,F) = \sum_{x\in E} \sum_{y \in F} (A\circ A) (x-y) = \sum_z (E\circ F) (z) \cdot (A\circ A) (z) \ge \a |E| |F| \,.
\end{equation}
    We shall call a set $A$ is {\it strongly connected}
    if inequality (\ref{f:pre_strong}) holds.
    As was noted above in many problems of combinatorial number theory
    the order of the number $h$ was $T_2 (A)/|A|^2$.
    We also put $\a = C \cdot T_2 (A)/|A|^2$,
    where $C>0$ is a constant.

\Def
\label{a_connectedness}
    Let $k\ge 2$ be a positive integer.
    An arbitrary nonempty finite set $A\subseteq G$ is called {\it strongly connected of degree $k$}
    if there is an absolute constant $C\in (0,1]$
    such that
    for any disjoint sets
    $E,F\subseteq A$, $E \bigsqcup F = A$, we have
    \begin{equation}\label{ineq:a_connectedness}
        \sum_x (E \circ F)(x) \cdot ( (A *_{k-2} A) \circ (A *_{k-2} A) ) (x) \ge C c_E c_F T_k (A) \,,
    \end{equation}
    where $c_E = |E| / |A|$, $c_F = |F| /|A|$.

    First of all let us show that any strongly connected set is a connected set.

\St
\label{st:strong_implies_weak}
{\it
    Let $p$ be a positive integer, and $k = 2^p$.
    Suppose that $A$ is a strongly connected of degree $k$ set
    such that (\ref{ineq:a_connectedness}) holds with some constant $C$.
    Then $A$ is connected of degree $k$ and inequality (\ref{ineq:b_connectedness})
    holds with $C/8$.
}
\\
\Proof
If the cardinality of $A$ is less than two then there is nothing to prove.
Let $|A| \ge 3$, $B$ be an arbitrary subset of $A$, and $\ov{B} = A\setminus B$.
Let also
\begin{equation}\label{}
    \sigma = \sum_x (B \circ \ov{B})(x) \cdot ( (A *_{k-2} A) \circ (A *_{k-2} A) ) (x) \,.
\end{equation}
Since $A$ is a strongly connected of degree $k$ it follows that
\begin{equation}\label{tmp:19:54_28}
    \sigma \ge C \frac{|B|}{|A|} \frac{|\ov{B}|}{|A|} \cdot T_k (A) \,.
\end{equation}
We have
\begin{equation}\label{}
    \sigma = \sum_x (B * A *_{k-2} A) (x) \cdot (\ov{B} * A *_{k-2} A) (x) \,.
\end{equation}
Using Lemma \ref{l:T_k(B,A)}, we get
$\sigma^{2k} \le T_k (B) T_k (\ov{B}) T^{2k-2}_k (A)$.
Combining the last inequality and (\ref{tmp:19:54_28}), we obtain
\begin{equation}\label{}
    T_k (B) T_k (\ov{B}) \ge C^{2k} \frac{|B|^{2k}}{|A|^{2k}} \cdot \frac{|\ov{B}|^{2k}}{|A|^{2k}} T^2_k (A) \,.
\end{equation}
If $|B| \le |A|/2$ then $|\ov{B}| \ge |A|/2$.
Using this lower bound for $|\ov{B}|$,
we get
\begin{equation}\label{tmp:intermediate}
    T_k (B) \ge \left( \frac{C}{2} \right)^{2k} \left( \frac{|B|}{|A|} \right)^{2k} T_k (A)
\end{equation}
and the statement is proved.
Suppose that $|B| > |A|/2$.
Then let $B_1$ be an arbitrary subset of $B$ of the cardinality $[|A|/2]$.
Clearly, $|B| \le 4 |B_1|$.
By (\ref{tmp:intermediate}), we have
$$
    T_k (B) \ge T_k (B_1)
        \ge
            \left( \frac{C}{2} \right)^{2k} \left( \frac{|B_1|}{|A|} \right)^{2k} T_k (A)
                \ge
                    \left( \frac{C}{8} \right)^{2k} \left( \frac{|B|}{|A|} \right)^{2k} T_k (A) \,.
$$
This completes the proof.

Thus any strongly connected set is connected.
In particular, Proposition \ref{pr:small_basis} is true for an arbitrary strongly connected set
and therefore any strongly connected set is contained in $\Span \L$
for some $\L$ with small cardinality.
Apparently,
it was S.V. Konyagin (see \cite{Konyagin_sum_prod}) who first proved that
an arbitrary strongly connected set
is economically contained in some special subgroup
(see also another variant of his statement in book \cite{Tao_Vu_book} p. 114, ex. 2.6.10).
We formulate his result in our terms
and give the proof for the sake of completeness.

\St
\label{st:Konyagin_connectedness}
{\it
    Let $k \ge 2$ be a positive integer.
    Let also $A\subseteq G$ be a strongly connected of degree $k$ set
    such that (\ref{ineq:a_connectedness}) holds with some constant $C$.
    Let
    $$
        S = \left\{ h \in G ~:~ ( (A *_{k-2} A) \circ (A *_{k-2} A) ) (x) \ge C \frac{T_k (A)}{|A|^2} \right\} \,.
    $$
    Then there is an element $a\in G$ such that $A \subseteq \langle S\rangle + a$,
    where $\langle S\rangle$ is the subgroup of $G$ generating by $S$.
}
\\
\Proof
Assume the converse.
Let $H = \langle S\rangle$ and let $A_1, \dots, A_r \subseteq A$
be intersections of $A$ with cosets of $H$.
If there are two nonempty intersections of cosets of $H$ with $A$, say,
$A_i$ and $A_j$, $i<j$, $i,j \in \{ 1, \dots, r \}$
then put $E = \bigsqcup_{l=1}^i A_l$ and $F = A\setminus E$.
Clearly, $E$ and $F$ are nonempty sets.
Since for any $e\in E$ and $f\in F$, we have $e-f \notin H$,
and, consequently, $e-f \notin S$ it follows that
$$
    \sum_x (E \circ F)(x) \cdot ( (A *_{k-2} A) \circ (A *_{k-2} A) ) (x)
        \le
$$
$$
        \le
            \sum_{x\notin S} (E \circ F)(x) \cdot ( (A *_{k-2} A) \circ (A *_{k-2} A) ) (x)
                < C |E| |F| \cdot \frac{T_k (A)}{|A|^2}
$$
with contradiction.
This completes the proof.

We prove an analog of Theorem \ref{t:connected_subset}
for strongly connected sets.

Let $E,F \subseteq A$ be sets.
Denote by  $e(E,F)$
the quantity
$\sum_x (E \circ F)(x) \cdot ( (A *_{k-2} A) \circ (A *_{k-2} A) ) (x)$.
Clearly, $e(E_1 \bigsqcup E_2, F) = e(E_1,F) + e(E_2,F)$.
Suppose that $E \subseteq A$ is an arbitrary set.
By $c_E$ denote
the ratio $|E| / |A|$.

We need in the following technical definition of strongly connected of degree $k$ sets.

\Def
{
    Let $k\ge 2$ be a positive integer.
    An arbitrary nonempty finite set $A\subseteq G$ is called {\it $\beta$---strongly connected of degree $k$}
    if there is an absolute constant $C\in (0,1]$
    and a set $B \subseteq A$, $|B| \ge \beta |A|$
    such that
    for any disjoint sets
    $E,F\subseteq B$, $E \bigsqcup F = B$, we have
    \begin{equation}\label{ineq:a_connectedness+}
        \sum_x (E \circ F)(x) \cdot ( (A *_{k-2} A) \circ (A *_{k-2} A) ) (x) \ge C c_E c_F T_k (A) \,.
    \end{equation}
}

Our next statement can be proved likewise Statement \ref{st:strong_implies_weak}.

\St
\label{st:strong_implies_weak+}
{\it
    Let $p$ be a positive integer, $k = 2^p$, and $\beta \in [0,1]$ be a real number.
    Let $A$ be a strongly $\beta$---connected of degree $k$ set
    and (\ref{ineq:a_connectedness+}) is true with some $C$
    and some  $B\subseteq A$, $|B| \ge \beta |A|$.
    Then $B$ is connected of degree $k$ set and inequality (\ref{ineq:b_connectedness})
    holds with $C\beta^2/8$.
}

A graph--theoretical variant of Lemma \ref{l:partition_with_large_T_k}
below was proved in \cite{Ruzsa_Elekes}.

\Lemma
\label{l:partition_with_large_T_k}
{\it
    Let $k\ge 2$ be a positive integer,
    $\eps_1\in [0,1]$ be a real number,
    and let $A\subseteq G$ be a finite set.
    Then there exists a partition of $A$ into disjoint sets $A_1,\dots,A_l$ such that \\
    $1)~$ For all $i,j\in \{ 1,\dots,l \}$, $i\neq j$, we have $e(A_i,A_j) \le \eps_1 c_{A_i} c_{A_j} T_k (A)$. \\
    $2)~$ For any $i\in \{ 1,\dots,l \}$ the set $A_i$ has the following property :
          for any disjoint sets $E,F\subseteq A_i$, $E \bigsqcup F = A_i$,
          we have $e(E,F) \ge \eps_1 c_E c_F T_k (A)$.

    Besides, the following inequality holds \\
    $3)~$ $\sum_{i=1}^l T_k (A_i) \ge T_k (A) \cdot (1- (2k-1) \eps_1)$.
}
\\
\Proof
Consider all partitions of $A$ into disjoint subsets $A_1, \dots, A_s$,
where $s$ is an arbitrary positive integer.
Select one for which the sum
\begin{equation}\label{f:->min}
    \sigma (A_1, \dots, A_s) = \sum_{1\ge i < j \le k} (~ e(A_i,A_j) - \eps_1 c_{A_i} c_{A_j} T_k (A) ~)
\end{equation}
is minimal.
By minimality of this partition, say $\{ A_1, \dots, A_l\}$,
we have $2)$.

Let us prove $1)$.
Suppose that for some $i,j\in \{ 1,\dots,l \}$, $i\neq j$ the following holds
$e(A_i,A_j) > \eps_1 c_{A_i} c_{A_j} T_k (A)$.
Constructing the new partition $\mathcal{P}$ of the set $A$,
$\mathcal{P} = \{ A_r \}_{r\neq i,j} \bigsqcup (A_i \bigsqcup A_j)$
and using the last inequality, we get
$$
    \sigma (\mathcal{P}) = \sigma (A_1, \dots, A_s)
                                - ( e(A_i,A_j) - \eps_1 c_{A_i} c_{A_j} T_k (A) )
                                    < \sigma (A_1, \dots, A_s) \,.
$$
with contradiction.

Prove that $1)$ implies $3)$.
Indeed,
$$
    T_k (A)
        =
            \sum_x (A *_{k-1} A)^2 (x)
                =
                    \sum_{i,j=1}^l \sum_x (A_i * A *_{k-2} A) (x) \cdot (A_j * A *_{k-2} A) (x)
                        =
$$
$$
                        =
                            \sum_{i=1}^l \sum_x (A_i * A *_{k-2} A) (x) \cdot (A_i * A *_{k-2} A) (x)
                                +
                                    \sum_{i,j=1,\, j\neq i}^l \sum_x (A_i * A *_{k-2} A) (x) \cdot (A_j * A *_{k-2} A) (x)
$$
$$
        =
            \sum_{i=1}^l \sum_x (A_i * A *_{k-2} A) (x) \cdot (A_i * A *_{k-2} A) (x)
                                +
                                    \sum_{i,j=1,\, j\neq i}^l
                                        \sum_x (A_i \circ A_j)(x) \cdot ( (A *_{k-2} A) \circ (A *_{k-2} A) ) (x)
$$
$$
        \le
                                                \sum_{i=1}^l \sum_x (A_i * A *_{k-2} A) (x) \cdot (A_i * A *_{k-2} A) (x)
                                                    +
                                                        \eps_1 \sum_{i,j=1}^l c_{A_i} c_{A_j} T_k (A)
        \le
$$
$$
        \le
            \sum_{i=1}^l \sum_x (A_i * A *_{k-2} A) (x) \cdot (A_i * A *_{k-2} A) (x)
                +
                    \eps_1 T_k (A) \,.
$$
Hence
$\sum_{i=1}^l \sum_x (A_i * A *_{k-2} A) (x) \cdot (A_i * A *_{k-2} A) (x) \ge (1-\eps_1) T_k (A)$.
Similarly,
$$
    \sum_{i=1}^l \sum_{j=1}^l \sum_x (A_i * A_j * A *_{k-3} A) (x) \cdot (A_i * A *_{k-2} A) (x)
        \le
$$
$$
        \le
            \sum_{i=1}^l \sum_x (A_i * A_i * A *_{k-3} A) (x) \cdot (A_i * A *_{k-2} A) (x)
                +
$$
$$
                +
                    \sum_{i=1}^l \sum_{j=1,\, j\neq i}^l \sum_x (A * A_j * A *_{k-3} A) (x) \cdot (A_i * A *_{k-2} A) (x)
$$
$$
        \le
            \sum_{i=1}^l \sum_x (A_i * A_i * A *_{k-3} A) (x) \cdot (A_i * A *_{k-2} A) (x)
                +
$$
$$
                +
                    \sum_{i=1}^l \sum_{j=1,\, j\neq i}^l
                        \sum_x (A_i \circ A_j)(x) \cdot ( (A *_{k-2} A) \circ (A *_{k-2} A) ) (x)
        \le
$$
$$
        \le
            \sum_{i=1}^l \sum_x (A_i * A_i * A *_{k-3} A) (x) \cdot (A_i * A *_{k-2} A) (x)
                +
                    \eps_1 T_k (A) \,.
$$
And so on.
Finally, we obtain
$$
    \sum_{i=1}^l T_k (A_i)
        =
            \sum_{i=1}^l \sum_x (A_i *_{k-1} A_i) (x) \cdot (A_i *_{k-1} A_i) (x)
                \ge (1 - (2k-1) \eps_1) \cdot T_k (A) \,.
$$
This concludes the proof.

\Note
It is easy to see that the third property of the constructed partition implies that
there is $i_0 \in \{ 1, \dots, l \}$ such that
$|A_{i_0}| \ge (1-(2k-1) \eps_1) \cdot m^{\frac{\zeta_k (A) - 1}{2k-2}} \ge (1-(2k-1) \eps_1) \cdot m^{1/2}$.
In fact, we have
$$
    m^{\zeta_k (A)} (1-(2k-1) \eps_1) \le \sum_{i=1}^l T_k (A_i) \le ( \max_{i=1,\dots,l} |A_i| )^{2k-2} \sum_{i=1}^l |A_i|
    \le ( \max_{i=1,\dots,l} |A_i| )^{2k-2} m \,.
$$
This yields that if we put $\beta = (1-(2k-1) \eps_1) m^{-1/2}$ then  {\it any}
set $A\subseteq G$, $|A| = m$
is strongly $\beta$---connected of degree $k$
and inequality (\ref{ineq:a_connectedness+}) holds with any $C \le 1/ (2k-1)$.
Thus to obtain nontrivial results on the structure of $A$ one should
prove that $A$ is strongly $\beta$---connected for {\it large} $\beta$.

\Th
\label{t:strong_connected_subset}
{\it
    Let $A\subseteq G$ be a set.
    Let also $\eps, \beta \in (0,1)$ be real numbers, and $|A| \ge \eps / (2 \beta^2)$.
    Then there exists a partition of $A$ into disjoint sets $A_1, \dots, A_t$, $\Omega$ such that \\
    $1)~$ Any set $A_i$, $i=1,\dots,t$ is strongly $\beta$---connected of degree $2$
          and inequality (\ref{ineq:a_connectedness+}) holds with any $C \le \eps \log (1/\beta) / (3 \log (2|A|/\eps))$. \\
    $2)~$ $\sum_{i=1}^t T_2 (A_i) \ge (1-\eps) \cdot T_2 (A)$.
}
\\
\Proof
    Let $m = |A|$, $s_0 = \log (2m/\eps) / (2 \log (1/\beta)) \ge 1$,
    and $\eps' = \eps / (6 s_0)$.
    Let $C\le \eps'$ be a real number.
    The proof of Theorem \ref{t:strong_connected_subset} is a sort algorithm.
    If $A$ is strongly $\beta$---connected of degree $2$
    and (\ref{ineq:a_connectedness+}) is true with the constant $C$
    then there is nothing to prove.
    Suppose that $A$ is not strongly
    $\beta$---connected of degree  $2$ set (with the constant $C$).
    Using Lemma \ref{l:partition_with_large_T_k} with $\eps_1 = \eps'$, we get
    the partition $\mathcal{P}^{(1)}$ of $A$ into
    $A_1,\dots, A_l$ satisfy properties $1)$ --- $3)$ of the lemma.
    Since $A$ is not strongly $\beta$---connected of degree $2$
    it follows that for any $l\in \{ 1,\dots, l\}$, we have $|A_i| < \beta |A|$.
    Using the third property of the partition $\mathcal{P}^{(1)}$, we obtain
    $$
        \sum_{\mathcal{A} \in \mathcal{P}^{(1)}} T_2 (\mathcal{A}) = \sum_{i=1}^l T_2 (A_i) \ge (1-3\eps') T_2 (A) \,.
    $$
    Let $B^{(1)} = \{ A_i \mbox{ --- is not strongly } \beta-\mbox{connected of degree 2} \}$,
    and $G^{(1)}$ be the collection of all other sets of the partition $\mathcal{P}^{(1)}$.
    Let us construct a new partition of $A$.
    We do not change the sets $A_i$ from $G^{(1)}$.
    Further, for any $A_i$ belongs to $B^{(1)}$, we use Lemma \ref{l:partition_with_large_T_k}
    with $\eps_1 = \eps'$.
    We get a new partition of $A_i$ into subsets $A_{ij}$, $j\in \{ 1, \dots, l(i) \}$.
    So we construct a new partition $\mathcal{P}^{(2)}$ of the set $A$.
    For any $A_i \in B^{(1)}$ the following holds
    $\sum_{j=1}^{l(i)} T_2 (A_{ij}) \ge (1-3\eps') T_2 (A_i)$.
    Hence
    \begin{equation}\label{f:odin}
        \sum_{\mathcal{A} \in \mathcal{P}^{(2)}} T_2 (\mathcal{A}) \ge (1-3\eps')^2 \cdot T_2 (A) \,.
    \end{equation}
    Let
    $B^{(2)} = \{ A_{ij} \mbox{ --- is not strongly } \beta-\mbox{connected of degree 2} \}$.
    For an arbitrary $A_{ij} \in B^{(2)}$, we use Lemma \ref{l:partition_with_large_T_k}
    with $\eps_1 = \eps'$.
    We get a new partitions of the sets $A_{ij}$ into disjoint subsets $A_{ijr}$.
    And so on.
    At $s$--th step of the algorithm, we construct the partition $\mathcal{P}^{(s)}$
    such that
    \begin{equation}\label{f:posle_odin}
        \sum_{\mathcal{A} \in \mathcal{P}^{(s)}} T_2 (\mathcal{A}) \ge (1-3\eps')^s \cdot T_2 (A) \ge (1-3\eps' s) \cdot T_2 (A) \,.
    \end{equation}
    It is easy to see that if for some $s\le s_0$ the following holds
    \begin{equation}\label{f:tri}
        \sum_{\mathcal{A} \in \mathcal{P}^{(s)} \setminus B^{(s)}} T_2 (\mathcal{A}) \ge (1-\eps) \cdot T_2 (A) \,,
    \end{equation}
    then we are done.
    Indeed, just put $\Omega = \bigsqcup_{\mathcal{A} \in B^{(s)}} \mathcal{A}$.
    Suppose that for all $s\le s_0$ inequality (\ref{f:tri}) does not hold.
    Using inequality $s\le s_0$ and (\ref{f:posle_odin}), we get
    $
     \sum_{\mathcal{A} \in \mathcal{P}^{(s)}} T_2 (\mathcal{A}) \ge (1-\eps /2) \cdot T_2 (A)
    $.
    Hence
    \begin{equation}\label{f:chetire}
        \sum_{\mathcal{A} \in B^{(s)}} T_2 (\mathcal{A}) \ge \frac{\eps}{2} \cdot T_2 (A) \ge \frac{\eps m^2}{2} \,.
    \end{equation}
    For any $\mathcal{A} \in B^{(s)}$, we have $|\mathcal{A}| < \beta^s m$.
    Whence
    \begin{equation}\label{f:dva}
        \sum_{\mathcal{A} \in B^{(s)}} T_2 (\mathcal{A})
            <
                (\beta^s m)^2 \sum_{\mathcal{A} \in \mathcal{P}^{(s)}} |\mathcal{A}|
                    =
                        \beta^{2s} m^3 \,.
    \end{equation}
    If $s=s_0$ then the last inequality contradicts (\ref{f:chetire}).
    So for some $s< s_0$ inequality (\ref{f:tri}) holds.
    This completes the proof.

There is a difference between Theorem \ref{t:connected_subset}
and Theorem \ref{t:strong_connected_subset}.
In Theorem \ref{t:strong_connected_subset}
we prove that there exists a {\it partition} of $A$ into
strongly $\beta$---connected components and some
exceptional set $\Omega$
while Theorem \ref{t:connected_subset}
states that there is {\it one} connected subset of $A$.
Besides, Theorem \ref{t:strong_connected_subset}
implies that the remaining  set $\Omega$
has small $T_2 (\Omega)$.
Indeed, by the property $2)$, we have
$\sum_{i=1}^t T_2 (A_i) \ge (1-\eps) \cdot T_2 (A)$,
whence
$T_2 (\Omega) \le \eps T_2 (A)$.

\end{document}